\documentclass[12pt]{amsart}

\usepackage{fullpage}
\usepackage{hyperref}
\usepackage{amssymb,amsfonts,amsmath}

\usepackage{amssymb}
\usepackage{epstopdf}

\usepackage[all,cmtip]{xy}

\usepackage[T1]{fontenc}
\usepackage[utf8]{inputenc}
\usepackage{mathtools} 

\usepackage[bbgreekl]{mathbbol}

\usepackage{tikz}
\usepackage{color}
\usepackage{graphicx}
\usepackage{subfigure}
\graphicspath{{./Figure/}}

\newtheorem*{main theorem}{Main Theorem}

\theoremstyle{definition}
\newtheorem{construction}{Construction}
\newtheorem{comp}{Computation}
\newtheorem{obs}{Observation}

\theoremstyle{plain}
\newtheorem{thm}{Theorem}[subsection]

\newtheorem{lem}[thm]{Lemma}

\newtheorem{prop}[thm]{Proposition}
\theoremstyle{definition}
\newtheorem{defn}[thm]{Definition}
\newtheorem{rem}[thm]{Remark}
\newtheorem{prop/def}[thm]{Proposition/Definition}
\newtheorem{question}{Question}

\newcommand{\half}{\frac{1}{2}}
\newcommand{\dual}{\mathcal{R}_{[f]}^G}
\newcommand{\tors}{\mathbb{F}_{2}T_M}
\renewcommand{\int}{\operatorname{int}}

\newcommand{\Id}{\operatorname{Id}}

\newcommand{\sx}{{\sf{x}}}
\newcommand{\sy}{{\sf{y}}}

\newcommand{\Z}{\mathbb{Z}}

\newcommand{\fq}{\operatorname{fq}}
\newcommand{\id}{\operatorname{id}}

\newcommand{\stb}{S^2\times S^2}
\title{stabilizations of $s$-cobordisms of dimension $5$}

\author[Jinzhou Huang]{Jinzhou Huang}
\email{2100010619@stu.pku.edu.cn}

\begin{document}

\maketitle

\begin{abstract}

It has long been known that the $s$-cobordism theorem fails for $5$-dimensional $s$-cobordisms. In this article we study how many times of "stabilizations" are needed to turn a $5$-dimensional $s$-cobordism to a product cobordism. The question is analogous to asking how many times of stabilizations are needed to turn an exotic pair of four manifolds into diffeomorphic ones. The main tools in this article are Gabai's $4$D light bulb theorem and its applications, and we also use a refinement of $4$D light bulb theorem by Freedman Quinn invariant.
\end{abstract}

\section{Introduction}\label{sec1}

In this article we will study the "stabilization" of a $5$-dimensional $s$-cobordism. First let us focus on $h$-cobordisms.  Wall 's result concerning the stabilizations of exotic $4$-manifolds has been well known to us \cite{WT}: Any two exotic simply-connected four manifolds become diffeomorphic after finite many times of stabilizations. Wall's strategy to prove this is to construct an $h$-cobordism between the two exotic four manifolds. Given Wall's result, people have been wondering how many times of stabilizations are needed to turn two exotic closed $4$-manifolds into diffeomorphic ones. However, we still don't know if one time is enough for this sake. Similar to the closed case, Gompf proved in \cite{GPF} that any two homeomorphic, compact, smooth orientable  four manifolds become diffeomrophic after sufficiently many stabilizations. Lately Sungkyung Kang found example that there exist homeomorphic smooth contractible $4$-manifolds $W_1$, $W_2$, with diffeomorphic boundaries, such that $W_1\#\stb$ and $W_2\#\stb$ are not diffeomorphic\cite[Corollary 1.2]{KANG}. We will not touch these questions in this article. Similar to the existence of exotic $4$-manifolds, Donaldson proved in \cite{SKD} that five dimensional $h$-cobordism is not necessarily a product cobordism , so we can ask a similar question regarding $h$-cobordisms: 

\begin{question}\label{que}
	
	\

How many times of "stabilizations" do we need to turn an $h$-cobordism into a trivial one? Here a stabilization of an $h$-cobordism is to make connected sum with a copy of $S^2\times S^2$ at each level of the original cobordism(if we are given a compatible Morse function). For more detailed definition the reader can refer to Definition \ref{stab}. 
\end{question}
As a related topic we introduce the complexity of an $h$-cobordism: Given an $h$-cobordism $(W;M_0,M_1)$ between smooth simply-connected four manifolds, we define the complexity $C(W;M_0,M_1)$ following \cite{MO} by:
\[
\begin{split}
	C(W;M_0,M_1)=&\operatorname{min}\{C(f)|f\ \textrm{is a compatible Morse function for $W$ with only critical} \\&\textrm{points of index $2,3$}\}
\end{split}
	\]

Here $C(f)$ is the number of flow lines from index $3$ points to index $2$ points minus the number of index $2$ points.

For a partial answer to the Question \ref{que}, an exercise \cite[Page 107]{FQ} in Freedman and Quinn's book implies that finite times are enough, however the number implied in their book is the total number of intersections among the Whitney disks which pair all the extra intersections of descending spheres and ascending spheres, which is no less than $C(W;M_0,M_1)$.

There is another formulation of the complexity of an $h$-cobordism, that is the minimal number of pairs of $2,3$ critical points for a compatible Morse function $f$ with only $2,3$ index points, which we denote by $C'(W,M_0,M_1)$. The main result of this article is that $C'(W;M_0,M_1)$ times of stabilizations are enough for Question \ref{que}, that is:
\begin{thm}\label{main}
	
	\
	
	Given an oriented $5$-dimensional s-cobordism $W=W^5$ between oriented, closed $4$-dimensional manifolds $M_0$ and $M_1$ and a Morse function $f:W\to\mathbb{R}$ satisfying that $f$ has only critical points of index $2$, $3$.
	
	Assume $f$ has $k$ pairs of critical points of index $2$, $3$, then $W$ becomes a product cobordism after $k$ times of stabilizations. 
\end{thm}

The key ingredient of the proof of Theorem \ref{main} is the following theorem \ref{tech0}. The key strategy of Smale's proof of $h$-cobordism theorem in higher dimensions is to isotope the descending spheres to a position where it intersects the corresponding ascending spheres in exactly one point via Whitney moves. Although we can not do the same thing in an $h$-cobordism between four manifolds, we prove that such a "good position" still exists(under certain assumptions), and many topoloical invariants obstructing an isotopy between the original sphere and the good position vanish. Moreover, the good position is unique up to isotopy.

Here we make a remark about the above definitions of complexity of an $h$-cobordism: As mentioned by John W. Morgan and Zolt\'{a}n Szab\'{o} in \cite{MO}, there is a family of $h$-cobordisms with arbitrary large $C$ while with $C'\equiv1$. As an result of the main theorem \ref{main}, an $h$-cobordism between simply-connected diffeomorphic $4$ manifolds with indefinite intersection forms becomes a product cobordism after a single stabilization (we do not really need \ref{main} to prove this actually, see Remark \ref{app}). This can be seen as a partial answer to the "one is enough" question concerning $h$-cobordisms rather than exotic $4$ manifolds. 

\begin{thm}\label{tech0}
	
	\
	
	Given an oriented $4$-dimensional manifold $M$ and smoothly embedded spheres $S_0^i$ and $S_1^i(i=1,\hdots,n)$ satisfying for every $i,j=1,\cdots,n$:
	\begin{enumerate}
		\item[C1]$\lambda_2(S_0^i,S_1^j)=\pm\delta_{i,j}[1]$, where $1$ represents the trivial element in $\pi_1$. 
		
		\item[C2]$S_0^i$ has trivial normal bundle.
		
	    \item[C3]$S_0^{i}\cap S_0^{j}=\varnothing$, $S_1^{i}\cap S_1^{j}=\varnothing$ for $i\ne j$. 
	    
		\item[C4]$\pi_1(M\backslash \cup_{i=1}^{n} S_0^i)\to\pi_1(M)$ and $\pi_1(M\backslash \cup_{i=1}^{n} S_1^i)\to\pi_1(M)$ are both injective. 
	\end{enumerate}
	Then there is another family of spheres $S_3^i(i=1,...,n)$ satisfying for $i,j=1,\cdots,n$:
	\begin{enumerate}
		\item[R1]$S_3^i(i=1,\hdots,n)$ are disjointly embedded.
		
		\item[R2]$S_3^i$ intersects $S_0^j$ transversely and $\#S_3^i\cap S_0^j=\delta_{ij}$. 
		
		\item[R3]$S_3^i$ is homotopic to $S_1^i$.
		
		\item[R4]The Freedman-Quinn invariant $\fq(S_1^i, S_3^i)=0$ for each $i=1,\hdots,n$. 
	\end{enumerate}
	Moreover the resulted family $S_3$ is unique up to isotopy, that is if $S_3'^i(i=1,\hdots,n)$ is another family of spheres satisfying the above four conditions, then there exists an isotopy from $\cup_{i=1}^{n}S_3^i$ to $\cup_{i=1}^{n}S_3'^i$.
\end{thm}

The main ingredient of the proof of above two theorems is Gabai's $4$-D light bulb theorem \cite[Theorem 1.2]{DG}, which is enough for us to prove the main theorem in the simply-connected situation. However, when the fundamental group is non-trivial, we must use a refinement of Gabai's $4$-D light bulb theorem by Freedman Quinn invariant in Rob Schneiderman and Peter Teichner's article\cite{SRP}. Dave Auckly, Hee Jung Kim, Paul Melvin, Daniel Ruberman and  Hannah Schwartz proved in \cite{AD} that under certain assumptions two homotopic embedded spheres in a four  manifold become isotopic in the once-stabilized manifold, we will also borrow a lot of techniques in their article.

All manifolds in this article are smooth.
\

{\bf Acknowledgment.} I would like to thank Professor Jianfeng Lin for telling me this problem, who had some very inspirational conversations with me, offered many valuable advises on this paper and provided a lot of useful materials for me to study and refer to.

\

\section{Preliminaries}
\subsection{Tools In $4$-Dimensional Topology}\label{tool}

\

In this subsection $M$ denotes a $2n$-dimensional smooth oriented manifold(possibly with boundaries) with a chosen basepoint, and $A$, $B$ are \emph{properly immersed} oriented \emph{simply-connected} manifolds(possibly with boundaries). We further assume all intersections of $A$ and $B$ are transverse and are double points. The main source of this subsection is \cite[Chapter 1]{FQ}.

We have already been familiar with algebraic numbers with values in $\mathbb{Z}$. In case where $A$, $B$, $M$ are all closed manifolds, the algebraic intersection number of $A$ and $B$ in $M$ is: \[\langle[A]\smile[B], [M]\rangle=\sum_{p\in A\cap B}\epsilon(p) \] where $\epsilon(p)\in\{1,-1\}$ is the local intersection number, and $[M]$ is the fundamental class.

When the fundamental group of $M$ is non-trivial, we get more information about the intersections of two properly immersed surfaces.

\begin{defn}[intersection number with values in $\mathbb{Z}\pi_1(M)$]\label{intersection}
	
	\

 We pick two arcs connecting the base point $z\in M$ to $w\in A$ and $w'\in B$ respectively, and we will call them \emph{whiskers}.

We define the intersection number $\lambda_n(A,B)\in\mathbb{Z}\pi_1(M)$ by adding up the elements $\sigma(p)g(p)$ for each intersection point $p$, where $\sigma(p)\in\{1,-1\}$ and $g(p)\in\pi_1(M)$. $\sigma(p)$ is the local intersection number at $p$. The element $g(p)$ is represented by a loop going along first the whisker from $z$ to $w$ , then a path in $A$(crossing no self-intersection points of $A$) from $w$ to $p$, then a path in $B$ from $p$ to $w'$(crossing no self-intersection points of $B$) and finally the whisker from $w'$ to $z$. Since $A$ and $B$ are simply-connected, $g(p)\in\pi_1(M)$ does not depend on the choices of the paths in $A$ and $B$ .

We can also define self-intersection number. Now for a self intersection $p$ of $A$, $g(p)$ is represented by the juxtaposition of the image of two paths in the domain of $A$ and $B$, where the first goes from $w$ to one of the preimage of $p$, and the second goes from another preimage to $w$, and we make $g(p)$ to be a loop based at $z$ via the whisker for $A$. In this case we can not order the two sheets at an intersection in a canonical way, so the self-intersection number $\mu_n(A)$ defined in the above procedure takes values in $\mathbb{Z}\pi_1/(g-(-1)^ng^{-1})$(Here $(g-(-1)^ng^{-1})$ denotes the \emph{abelian subgroup} generated by elements of the form $g-(-1)^ng^{-1}$ for $g\in \pi_1(M)$). We can further define $\tilde{\mu}_n(A)\in\mathbb{Z}\pi_1/(g-(-1)^ng^{-1},1)$ which is the equivalent class of $\mu_n(A)$
in $\mathbb{Z}\pi_1/(g-(-1)^ng^{-1},1)$.
\end{defn}

\begin{rem}

\

1. $(-1)^n$ in the definition of self-intersection number comes from the change of intersection number $\sigma$ when we switch the order of two sheets.

2. We can easily see that $\lambda_n$ and $\mu_n$ are invariants under regular homotopy(if we keep track of the whiskers): A regular homotopy can always be factored as isotopies, Whitney moves and finger moves \cite[\S 1.6]{FQ}, and both of the latter operations do not change $\lambda$ of $\mu$ since they eliminate or create a canceling pair $g$ and $-g$ respectively. $\tilde{\mu}_n$ is a invariant under homotopy, since a homotopy may result in a "cusp".

3. The intersection numbers defined above depend on the choices of whiskers. For example a different choice of whisker lead to a self-intersection number differing by a conjugation of an element in $\pi_1(M)$. 
\end{rem}

	\begin{figure}[!htb]
	\centering
	\subfigure[before tubing]{\resizebox*{8cm}{!}{\includegraphics{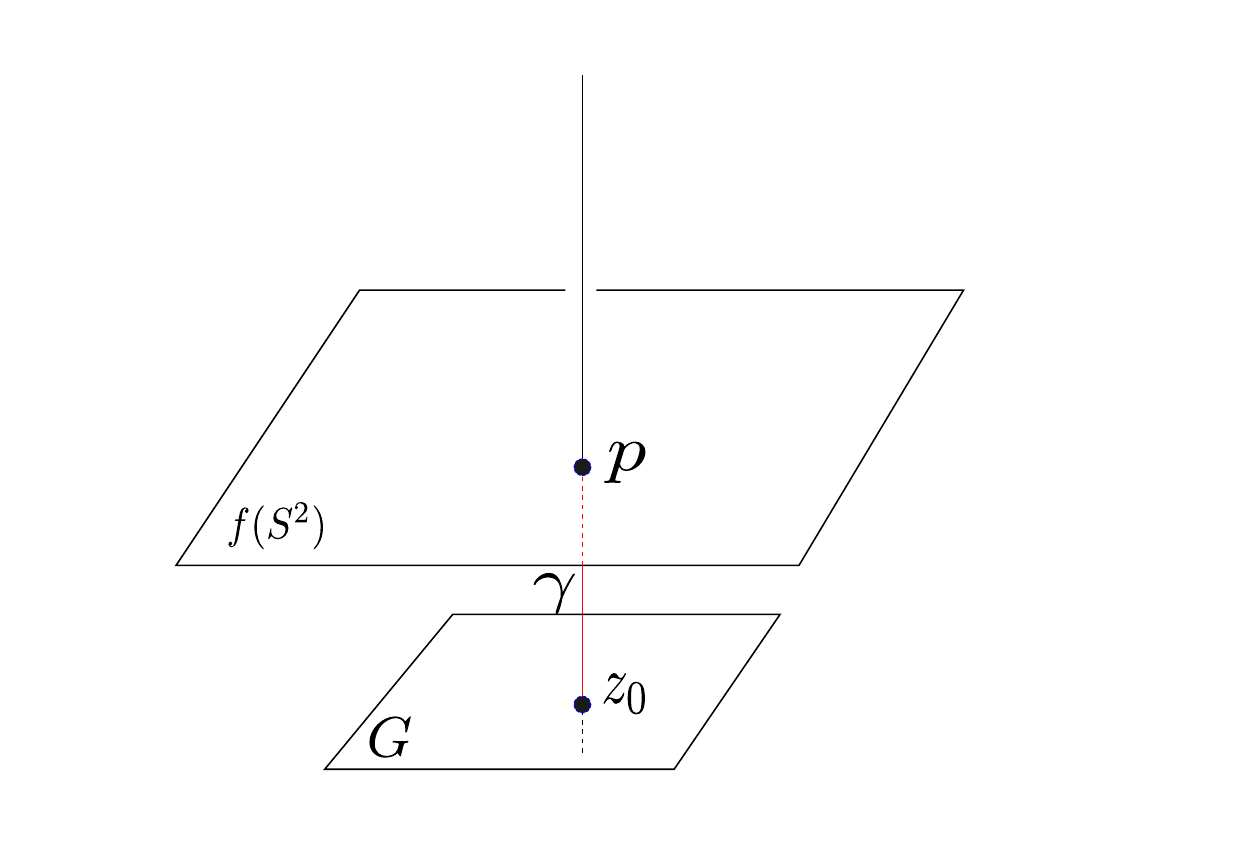}}}\hspace{2pt}
	\subfigure[after tubing]{\resizebox*{8cm}{!}{\includegraphics{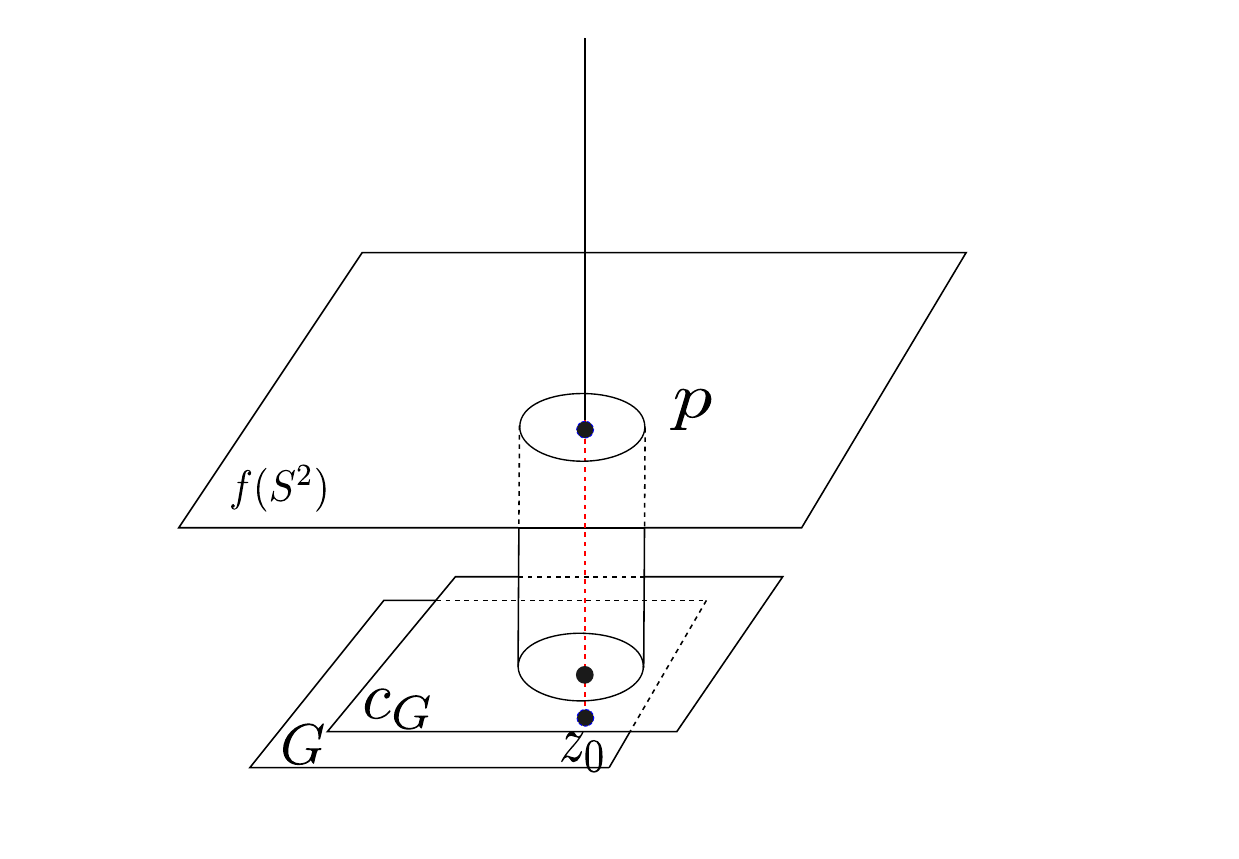}}}\hspace{2pt}
	\caption{}\label{figure}
\end{figure}

\begin{construction}[tubing]\label{tubing}
	
	\

	We introduce the well known technique of \emph{tubing} in $4$-dimensional topology. Assume $f:S^2\looparrowright M$ is an immersed sphere in $M$ with only transverse double self-intersections. $G$ is an embedded sphere which is a \emph{geometric dual} of $f$, that is $G$ has trivial normal bundle and intersects $f(S^2)$ transversely in exactly one point $z_0$. 
	
	For a self-intersection $p$ of $f$, by \emph{tubing off $p$ along a path $\gamma$ to a copy $c_G$ of $G$} we mean the operation indicated in Figure \ref{figure}. When we do multiple tubings the guiding paths may intersect, and in this situation we assign different radius for the tubes and take disjoint copies of $G$ so that the tubings do not produce extra intersections. For more details the reader can refer to \cite{NR} or \cite[\S 1.9]{FQ}.

\end{construction}

Here we mention a proposition about Whitehead torsion. Although Whitehead torsion has rich properties, we will only use an easy property concerning $s$-cobordisms here.

In this article an $h$-cobordism is an \emph{oriented} cobordism $(W;M_0,M_1)$ where the inclusions $M_0\hookrightarrow W$ and $M_1\hookrightarrow W$ are both homotopy equivalences. An $s$-cobordism is an $h$-cobordism with vanishing Whitehead torsion $\tau(W,M_0)$.

The following is indicated in the proof of \cite[Theorem 7.1D]{FQ} in Freedman and Quinn's book. For more information about Whitehead torsion the reader can refer to \cite{MIL2} or \cite{DP}.
\begin{prop}\label{White}
	
	\
	
If $(W;M_0,M_1)$ is an $h$-cobordism between $4$-dimensional manifolds $M_0$ and $M_1$, and $f:W\to \mathbb{R}$ is a compatible Morse function with only critical points of index $2$, $3$, and $\xi$ is a gradient-like vector field for $f$$($see \cite[Definition 3.1]{MIL}$)$. Then $\tau(W,M_0)=0$ is equivalent to that we can choose another compatible Morse function $f'($still with only critical points of index $2$ and $3$$)$ and a corresponding gradient-like vector field $\xi'$, and orientations and whiskers of the new descending and ascending spheres, so that the intersection numbers of new descending spheres $S_i$ and ascending spheres $T_j$ satisfy 
\[\lambda_2(S_i,T_j)=
\left\{
\begin{array}{ll}
0&(i\neq j)\\
1\in\pi_1&(i=j)
\end{array}\right.
\]
\end{prop}

\

\subsection{The Geometric Action and Freedman-Quinn Invariant}\label{geoaction}

\

In this subsection $M$ denotes an oriented smooth $4$-dimensional manifold without boundary. By the notation of \cite{SRP}, we denote 
\begin{displaymath}
\begin{split}
	\mathcal{R}_{[f]}^G=&\{\textrm{embedding $R: S^2\to M | R$ is homotopic to $f$ and has $G$ as a geometric dual}\} 
	\\ & /{\textrm{isotopies of $R$}}
	\end{split}
	\end{displaymath}

We define the $2$-torsions of $M$ by $T_M=\{g\in\pi_1(M)|g^2=1\neq g\}$, and $\mathbb{F}_2 T_M$ the $\mathbb{F}_2$-linear space generated by $T_M$.

In this subsection we will describe a transitive action of $\mathbb{F}_2 T_M$ on $\dual$ with known stabilizers $\mu_3(\pi_3(M))$, and we will describe the connections between Freedman-Quinn invariant and this action. Although the action described in \ref{action} is sufficient for usage in this article, we will still introduce Freedman-Quinn invariant for completion.

\

We now define Freedman-Quinn invariant following \cite[Definition 4.5]{SRP}.

\begin{defn}[Freedman-Quinn invariant]\label{FQdef}
	
	\
	
Given two embedded spheres $R$, $R'$ that are based homotopic in $M$, we choose a generic based homotopy $H: S^2\times I\to M\times\mathbb{R}$ from $R\times \{0\}$ to $R'\times \{0\}$. Denote $\hat{H}:S^2\times I\to M\times\mathbb{R}\times I$ the track of $H$, that is $\hat{H}(x,t)=(H(x,t), t)$. We define the Freedman Quinn invariant of $R$ and $R'$ to be $\fq(R,R')=[\mu_3(\hat{H})]\in\mathbb{F}_2 T_{M\times\mathbb{R}\times I}/\mu_3(\pi_3(M\times \mathbb{R}\times I))\cong\tors/\mu_3(\pi_3(M))$$($Since $S^2\times I$ is simply-connected, we can define the self-intersection number of $S^2\times I$ following Definition \ref{intersection}$)$.  
\end{defn}
\begin{rem}

\

\begin{enumerate}\label{remofFQ}
\item[1)]The definition of Freedman Quinn invariant is irreverent to geometric dual.

\item[2)]Intuitively Freedman-Quinn invariant counts how many times self-intersections occur in some time slice when we perform a generic homotopy from $R$ to $R'$.

\item[3)]Although generally $\mu_3$ takes values in $\mathbb{Z}\pi_1(M)/(g+g^{-1})$, we can prove that $\mu_3(\hat{H})$ actually takes value in $\tors<\mathbb{Z}\pi_1(M)/(g+g^{-1})$\cite[Lemma 4.1]{SRP}. 
\item[4)]For a pair of oriented four manifold $M,N$ where $N$ is simply-connected, $R,R'$ two homotopic spheres in $M$, then $\fq(R,R')$ is the same when calculated either in $M$ or $M\# N$(the connected sum is made away from $R,R'$): If we denote $B\subset M$ the ball where the connected sum is made, then by dimension argument the homotopy $H$ used to define Freedman Quinn invariant in Definition \ref{FQdef} can be made to be in  $(M\backslash B)\times\mathbb{R}$.
\end{enumerate}
\end{rem}

We can justify that $\fq$ does not depend on the choice of the based homotopy in Definition \ref{FQdef} by the following lemma from \cite[Lemma 4.4]{SRP}.

\begin{lem}
	\
	
If $J: S^2\times I \to M \times\mathbb{R}\times I$ is a generic track of a based self-homotopy of $R: S^2\to M\times\{0\}$, then $\mu_3(J)\in\tors$ lies in the image of the homomorphism $\mu_3: \pi_3(M\times{\mathbb{R}}\times I)\to\mathbb{F}_2 T_{M\times\mathbb{R}\times I}=\tors$.
\end{lem}

In presence of geometric dual, Freedman Quinn invariant can be used to distinguish isotopy classes. The following in from \cite[Corollary 1.6]{SRP}.

\begin{thm}\label{equal}
	Two elements are the same in $\dual$ if and only if their Freedman Quinn invariant vanishes.
\end{thm}
The following theorem is from \cite[\S 5.C]{SRP}.

\begin{thm}[The geometric action of $\tors$ on $\dual$]\label{transitive}
	
	\

	There is a transitive action of abelian group $\mathbb{F}_{2}T_M$ on $\mathcal{R}_{[f]}^G$ $($which we will denote by $\cdot$$)$, and $\mathcal{R}_{[f]}^G=\varnothing$ if and only if the reduced self-intersection number $\tilde{\mu}_2(f)$ does not vanish. The stabilizer of any $R\in \mathcal{R}_{[f]}^G$ is the subgroup $\mu_3(\pi_3 M)<\mathbb{F}_2 T_M$. If $R, R' : S^2\hookrightarrow M$ represent the same element in $\mathcal{R}_{[f]}^G$ and agree near $G$ then they are isotopic by an ambient isotopy supported away from $G$.
\end{thm}

\begin{rem}
	
	\
	
	We can derive Freedman-Quinn invariant from the above action: For $R$ and $R'$ in $\dual$, the Freedman-Quinn invariant
	$\fq(R,R')$ is $[t]\in\tors/\mu_3(\pi_3(M))$, where $t$ is any element in $\tors$ such that $t\cdot R=R'\in\dual$.\end{rem}

The following part of this subsection is devoted to describing explicit ways to calculate Freedman Quinn invariant and construct the action described in Theorem \ref{transitive}.
From now on we fix $f:S^2\looparrowright M$ with $\mu_2(f)=0$, and $f$ has a geometric dual $G$. We denote the self-intersection points of $f$ by $p_1,\hdots,p_{2n}$. Take $z$ the intersection of $f(S^2)$ with $G$ to be the basepoint, and the constant path at $z$ to be the whisker. The following definitions are from \cite{SRP}, 

\begin{defn}[choice of sheets]\label{sheet}
	
	\
	
	A \emph{choice of sheets} $\sx=\{x_1,\dots, x_{2n}\}\in$ for $f$ consists of choices $x_i\in f^{-1}(p_i)\subset S^2$ for each $i$(it can also be thought as a piece  $f(D)$ for $D$ a small disk neighborhood of $x_i$ in $S^2$). Each $x_i$ orients a loop $g_\sx (p_i)$ based at $p_i$ which is the image of the juxtaposition of two path in $S^2$, the first is the path from $z$ to $x_i$, the second is the path from $f^{-1}(p_i)$ to $z$, and via a path composed of first the chosen whisker and then \emph{a path} connecting $w$ and $p_i$ in $f(S^2)$.
\end{defn}

\begin{defn}\label{whi}
	\
	
  A \emph{clean collection of Whitney disks for $f$} is a set of Whitney disks $\{W_1,\dots, W_n\}$ which are framed, disjointly embedded, having interiors disjoint from $f$, disjoint from $G$, and pair all the self-intersection points of $f$.

We say a clean collection of Whitney disks $\{W_1,\dots, W_n\}$ are \emph{compatible} with the choice of sheets $\sx$ if for every Whitney disk $W_i$, there is \emph{a} Whitney arc $\alpha$ of $W_i$, such that $\alpha$ is the image of some arc in $S^2$ linking two elements in $\sx$.
  \end{defn}
  
  \begin{rem}\label{remwhi}
  	
  	\

Every clean collection of Whitney disks induces some choices of sheets compatible with the clean collection of Whitney disks(pick a component of the preimage of a Whitney arc for each Whitney disk and put its endpoints in $\sx$).

\end{rem}

\begin{obs}[The effect of tubing on the homotopy type of an immersed sphere]\label{eff}

\

At this point we make an observation about how tubing at a self-intersection $p$ to the geometric dual $G$ affects the homotopy class of the original immersed sphere $f$. We will also use this observation in Section \ref{sec2}.

\begin{figure}[!htb]
	\centering
	\subfigure[]{\resizebox*{8cm}{!}{\includegraphics{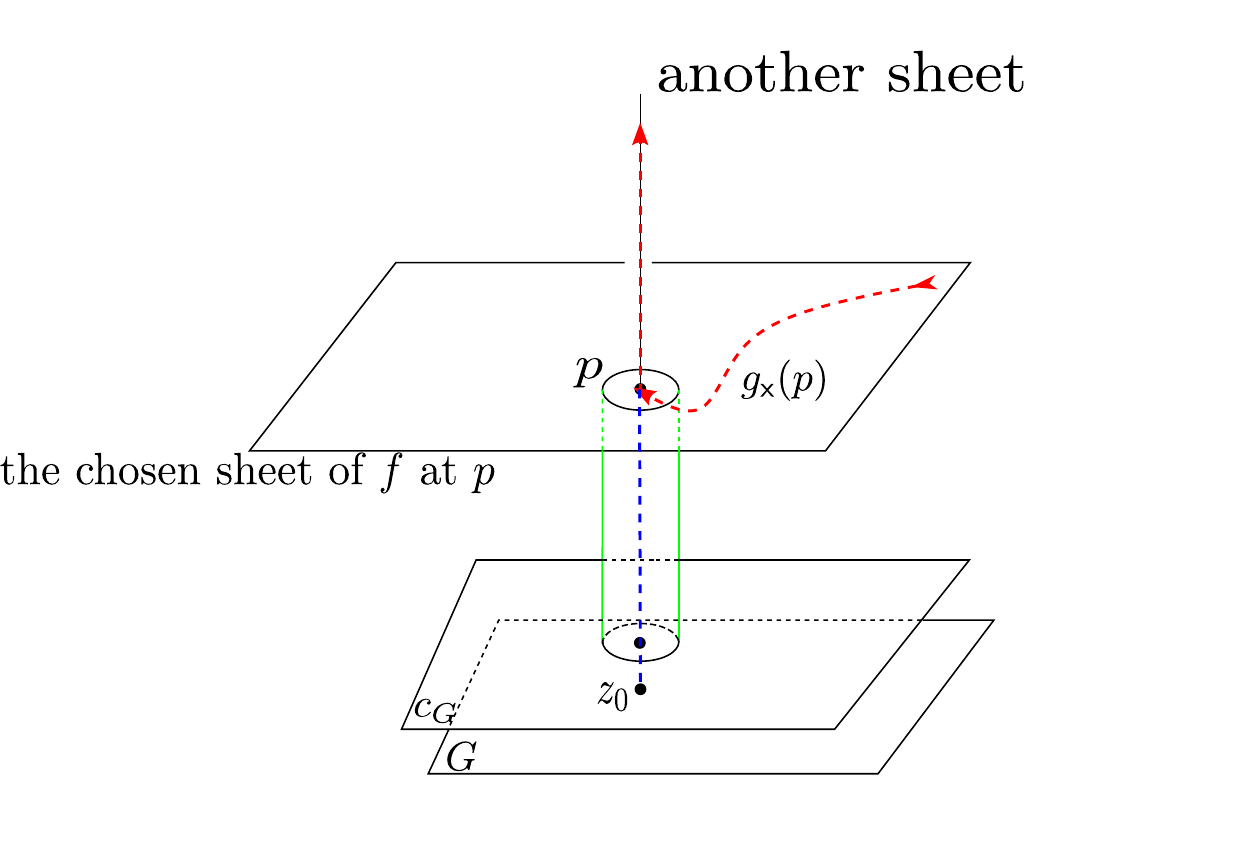}}}\hspace{2pt}
	\subfigure[moving along $g_{\sx}^{-1}(p)$]{\resizebox*{8cm}{!}{\includegraphics{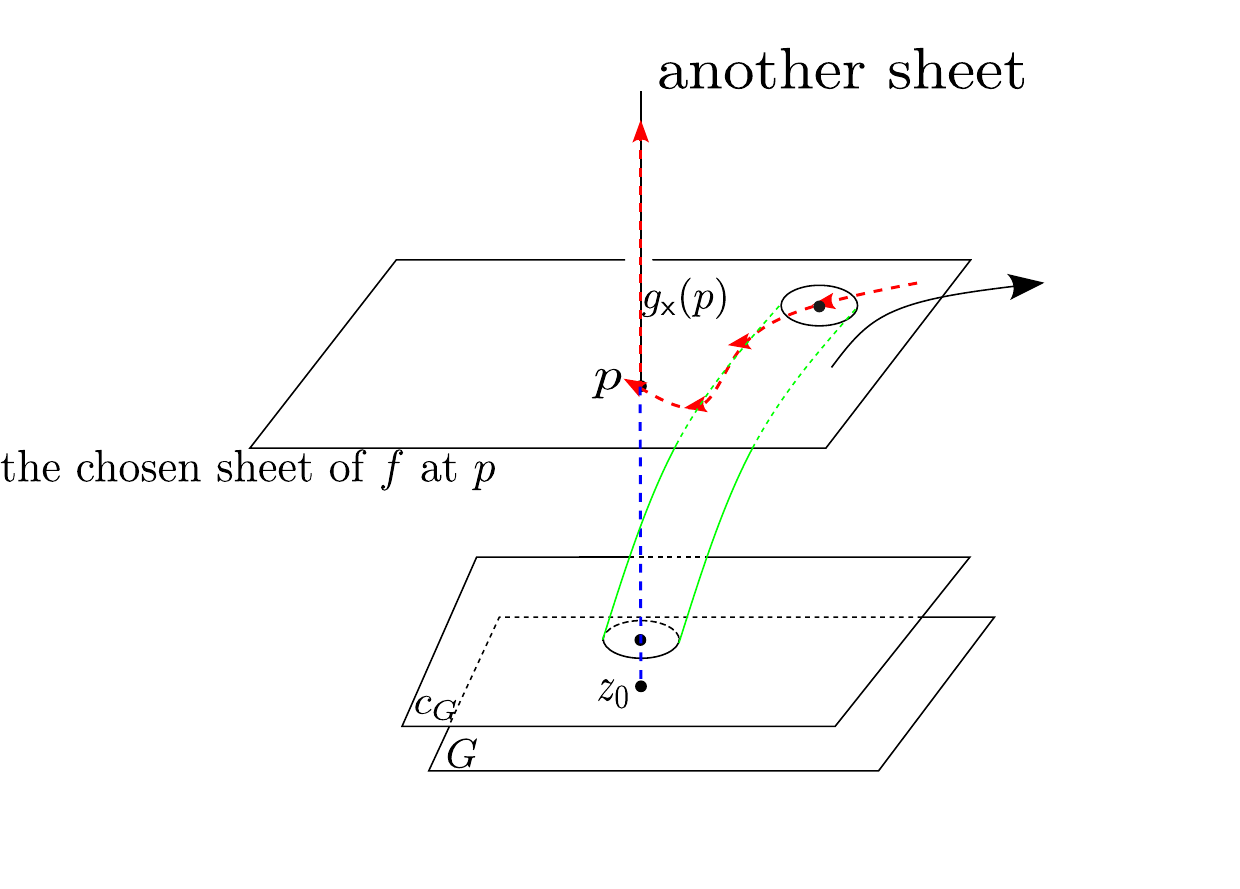}}}\hspace{2pt}
	\subfigure[]{\resizebox*{8cm}{!}{\includegraphics{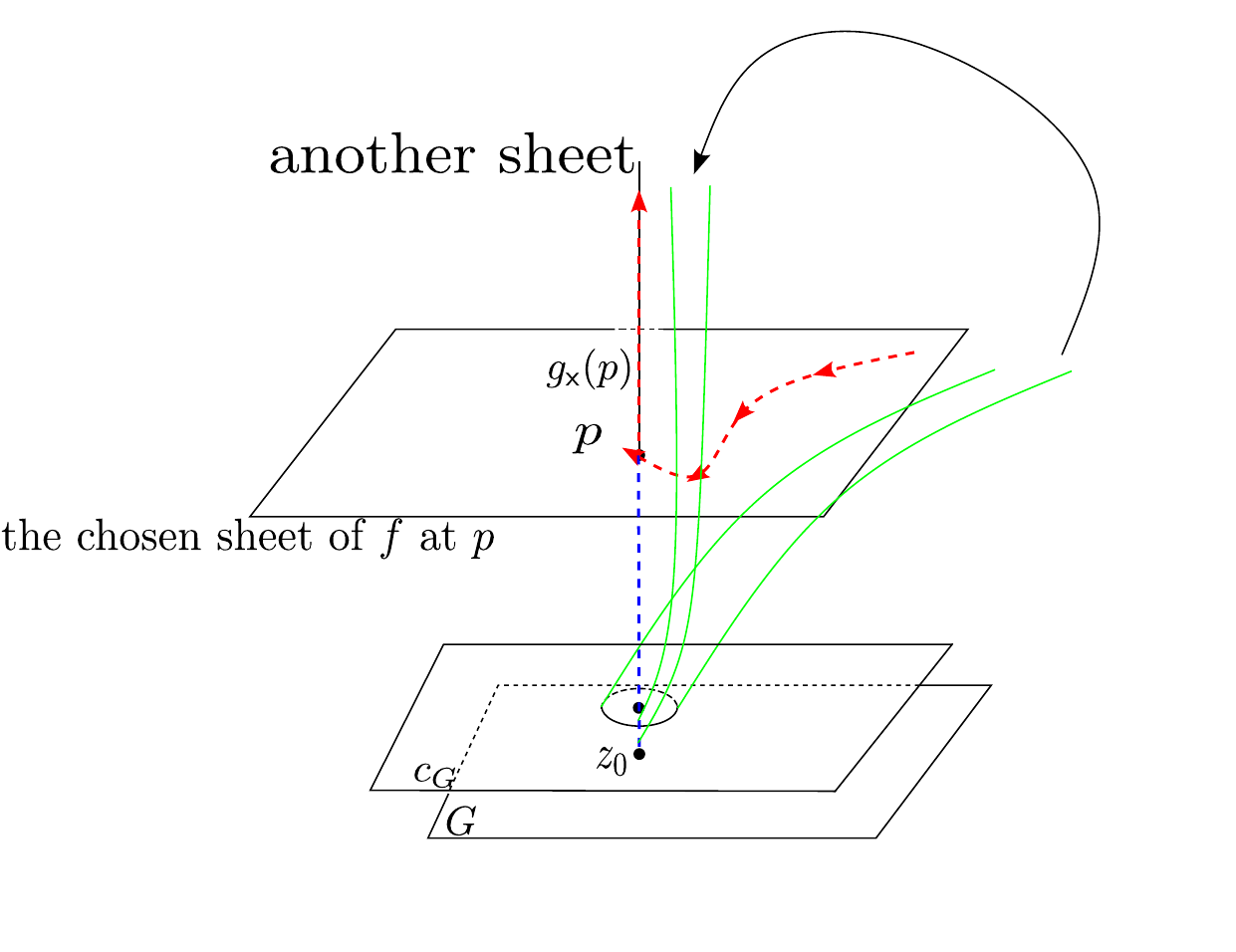}}}\hspace{2pt}
	\caption{}\label{proof}
\end{figure}
We orient $S^2$ and its geometric dual $G$ so that the intersection number at $z$ is $+1$. Recall the tubing operation in Construction \ref{tubing}.
\begin{prop}\label{propeff}
	
	\

For a choice of sheets $\sx$, we pick $x\in\sx$ and $p=f(x)$. We tube along a path which is the image of a path in $S^2$ from $f^{-1}(p)\backslash x$ to the preimage of $z_0$$($the blue dashed line in $($a$)$ of Figure \ref{proof}$)$. Then with the notation in Definition \ref{sheet}, the homotopy class of the result of tubing is $f+\epsilon(p)g_\sx(p)\cdot [G]\in \pi_2(M)$, where $\epsilon(p)$ is the intersection number at $p$.
\end{prop}

{\bf Proof:}
We indicate the proof in Figure \ref{proof}. First it is easy to see that the change of homology class is $\epsilon(p)[G]$ from (a) of Figure \ref{proof}. We move the tube along $g_{\sx}^{-1}(p)$ as indicated in (b) of Figure \ref{proof}.  Finally when we move the tube to the position in (c) of Figure \ref{proof} we easily see that the homotopy class is  $f+\epsilon(p)g_\sx(p)\cdot [G]$.\qed

\

If we do not want to change the homotopy class of $f$ after tubings, then our choice of sheets is required to satisfy
\begin{equation}\tag{$*$}\label{1}
0=\sum_{i=1}^{2n}\,\epsilon(p_i)\cdot g_\sx(p_i)\in\Z[\pi_1M]
\end{equation}
where $\epsilon(p_i)\in\{\pm 1\}$ is the local intersection number at $p_i$.
Such a choice of sheets exists since $0=\mu_2(f)\in\Z[\pi_1M]/(g-g^{-1})$.
\end{obs}

\

\cite{SRP} offers us a good way to compute and manipulate Freedman-Quinn invariant:

\

\begin{comp}[\bf Compute Freedman-Quinn invariant]

\

Choose a based regular homotopy $H$ from $R$ to $R'$ supported away from $G$, which begins with $n$ finger moves, following $n$ Whitney moves. Denote the stage after $n$ finger moves to be $h=H_{\half}$, then both $R$ and $R'$ can be obtained from $h$ by applying $n$ Whitney moves. The two Whitney moves induce two choices of sheets $\sx=\{x_1,...,x_{2n}\}$ and $\sx'=\{x_1',...,x_{2n}'\}$ as in Remark \ref{remwhi}. Each $x_i'\neq x_i$ gives an element $g_\sx (p_i)\in \pi_1(M)$, and we have the following proposition from \cite[\S 4.C]{SRP}.

\begin{prop}\label{compute}
	
	\

The Freedman Quinn invariant between $R$ and $R'$ can be calculated as 
\begin{equation}\tag{$**$}\label{2}
\fq(R,R')=\sum_{x_i\neq x_i'}g_\sx (p_i)\in\tors/{\mu_3(\pi_3(M))}
\end{equation}
\end{prop}
\

{\bf Proof:} We sketch the proof for the validity of the above computation. For more details the reader can refer to \cite[\S 4.C]{SRP}.

 $h$ is obtained from $R$ and $R'$ by $n$ finger moves. Following Definition \ref{FQdef}, we put the finger moves in $M\times\mathbb{R}$, and when a finger move in $M$ starts to generate self-intersections at a time, we use a bump function from $S^2$ to $\mathbb{R}$ to make the intersection to have different $\mathbb{R}$-coordinates. Finger moves(or their inverse Whitney moves) induce a compatible choice of sheets $\sx$, and we make the bump function positive in a neighborhood of elements of $\sx$, negative at elements in $f^{-1}(\{p_1,\hdots,p_{2n}\})\backslash\sx$. We denote the result by $h_\sx:S^2\hookrightarrow M\times\mathbb{R}$.

Now $R\times\{0\}$ and $R'\times\{0\}$ are isotopic in $M\times\mathbb{R}$ to $h_\sx$ and $h_{\sx'}$ for some $\sx$ and $\sx'$, respectively. When we homotope $h_\sx$ to $h_{\sx'}$ in $M\times\mathbb{R}$, we only need to scale the bump function for any element which is only in one of $\sx$ and $\sx'$ by $t$ from $1$ to $-1$, and this induces a self-intersection. Examine the elements in $\pi_1(M)$ which the self-intersection represents, and we complete the proof.\qed
\begin{figure}[!htb]
	\centering
	\subfigure[finger move along $t_i$]{\resizebox*{8cm}{!}{\includegraphics{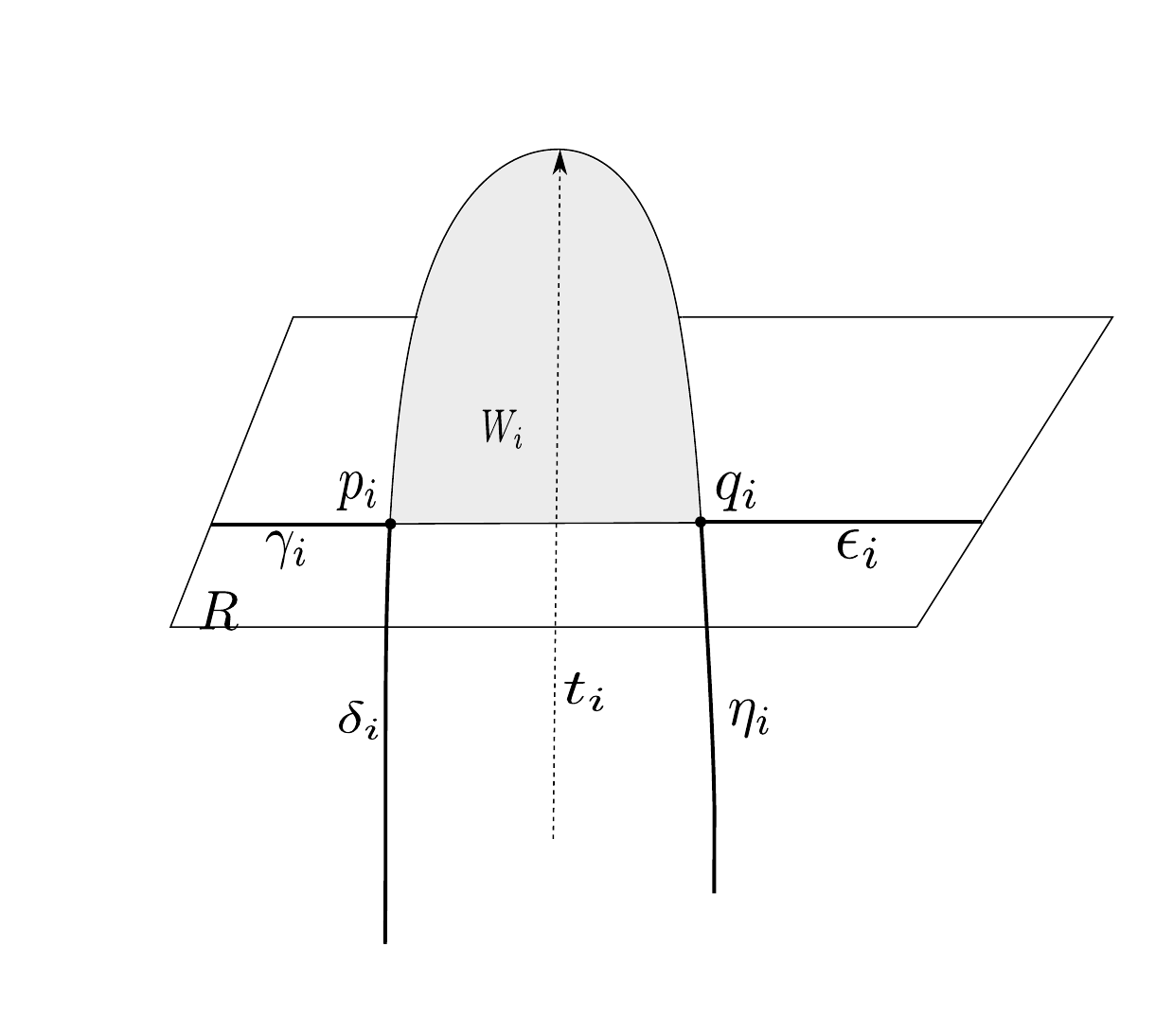}}}\hspace{2pt}
	\subfigure[tube along arcs indicated by $\delta_i$, $\eta_i$]{\resizebox*{8cm}{!}{\includegraphics{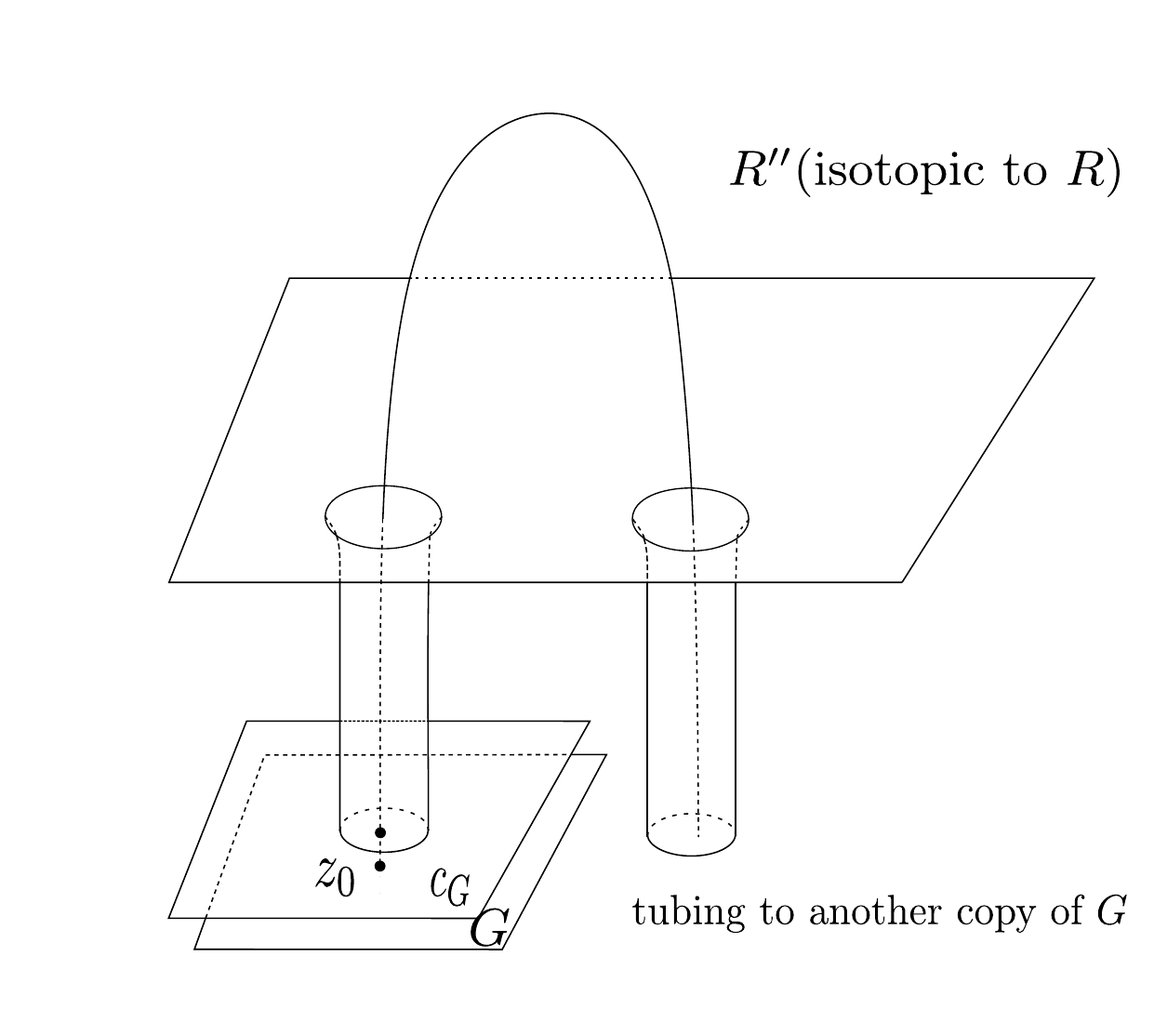}}}\hspace{2pt}
	\caption{}\label{figure2}
\end{figure}
\begin{figure}[!htb]
	\centering
	\subfigure[finger move along $t_i$]{\resizebox*{8cm}{!}{\includegraphics{the_disk.pdf}}}\hspace{2pt}
	\subfigure[tube along arcs indicated by $\delta_i$, $\epsilon_i$]{\resizebox*{8cm}{!}{\includegraphics{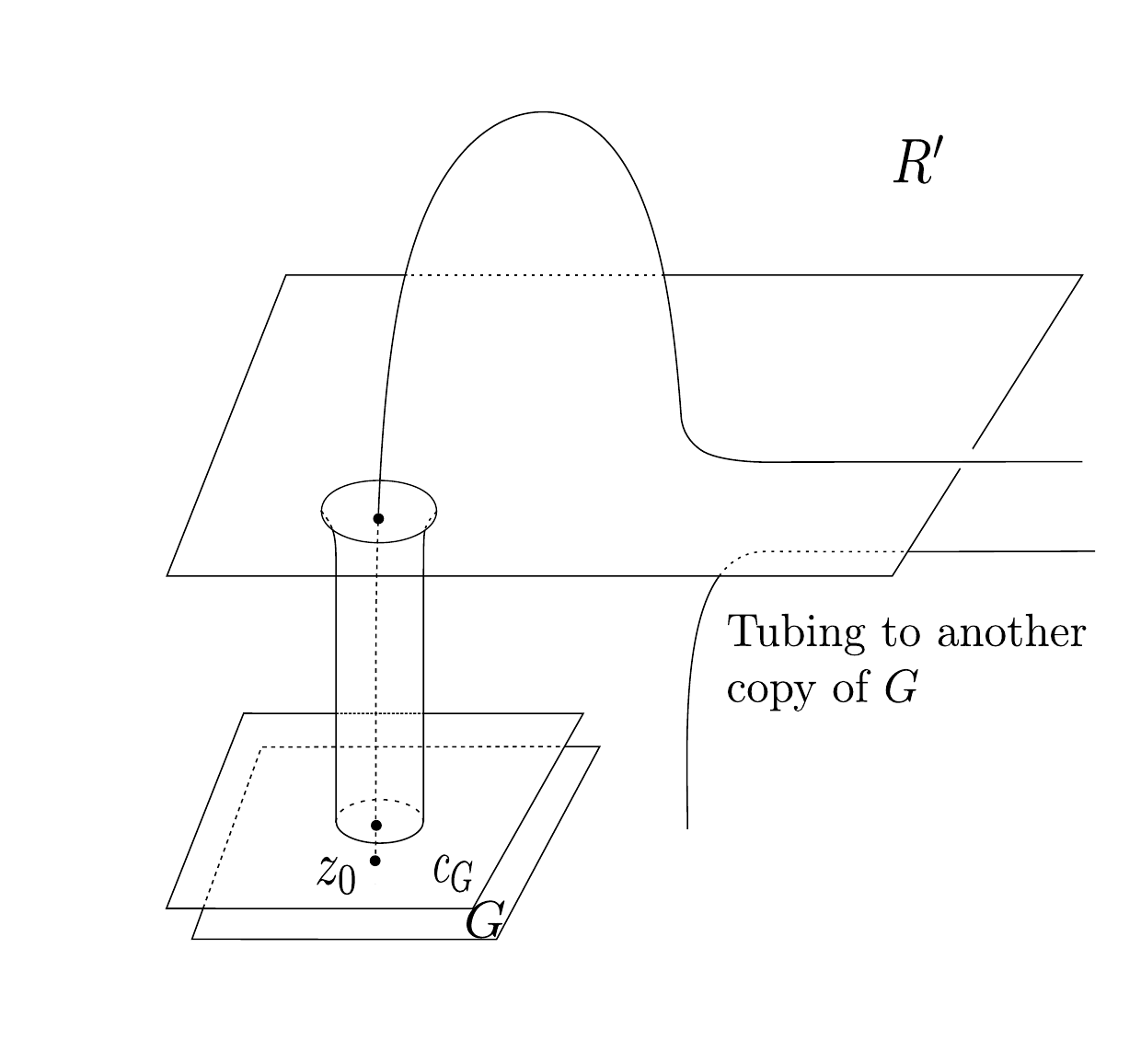}}}\hspace{2pt}
	\caption{}\label{figure3}
\end{figure}
\end{comp}
\

\begin{construction}[\bf $\tors$ action on $\dual$]\label{action}

\

Now we realize the transitive action mentioned in Theorem \ref{transitive} by a concrete geometric construction:

For $[R]\in\dual$ and $t\in\tors$, we define $R'$ as follows.
Write $t=t_1+\hdots+t_n$ with $t_i\in T_M$. For each $i=1,\hdots, n$, we perform finger moves along mutually disjoint arcs representing $t_i$ to create $n$ pairs of intersection points $p_i, q_i$ and Whitney disk $W_i$ as in the first picture of Figure \ref{figure2} and Figure \ref{figure3}, and we denote the resulted immersion by $h$. We tube along \emph{mutually disjoint} arcs $\delta_i$, $\eta_i$ to different copies of $G$ as in Figure \ref{figure3} for each pair of intersections $p_i,q_i$(since we are on an immersed sphere, we can always arrange the arcs to be disjoint), and we get an embedded sphere $R'$.(If we tube along arcs as in Figure \ref{figure2}, we get an embedded sphere isotopic to $R$).
\end{construction}
\

\begin{prop}
	
	\
	
	The sphere $R'$ constructed in the above procedure satisfies $\fq(R,R')=[t]\in\tors/\mu_3(\pi_3(M))$, so we can write $t\cdot R=R'$.
	\end{prop}
{\bf Proof:} We sketch the proof of $\fq(R,R')=[t]$. For more details the reader can refer to \cite[\S 5.A, \S 3.C]{SRP}.

For a general choice of sheets $\sy$ for $h$, for every $x\in\sy$ we tube along the arcs at each self-intersection point $p=h(x)$ to a copy of $G$, where the arc is the image of a path starting from $h^{-1}(p)\backslash\{x\}$ to the preimage of $z_0$ in $S^2$ as before (see the blue dashed line in (a) of Figure \ref{proof} for example). 
We use three facts:
\begin{enumerate}
\item[1)]The result of the above tubings is independent of the choices of the arcs used to guide the tubings. (see \cite[Lemma 3.5]{SRP}).
\item[2)]If a clean collection of Whitney disks is compatible with $\sy$, then there exists a choice of the arcs as above such that the result of tubing along these arcs is isotopic to the result of doing Whitney moves along such a clean collection of Whitney disks compatible with $\sy$.(see \cite[Lemma 3.3]{SRP})(the authors required in \cite[\S 3.C]{SRP} that the arcs used for tubings are disjointly embedded, which has been satisfied).
\item[3)]If $\sy$ satisfies (\ref{1}), then there exists a clean collection of Whitney disks compatible with $\sy$.(see \cite[Lemma 3.1]{SRP})
\end{enumerate}

We must remark that the presence of a geometric dual is vital for the above facts to be true.

Now $R$ is obtained from $h$ by Whitney moves which are inverses of the finger moves along $t_i$, and the Whitney moves induce a choice of sheets $\sx'$ as in Remark \ref{remwhi}.

We take a choice of sheets $\sx''$ by taking out from $\sx'$ the element in $h^{-1}(q_i)$ and putting back another element in $h^{-1}(q_i)$ for each $i$. $\sx''$ still satisfies (\ref{1}) since $t_i(i=1,\hdots, n)$ are $2$-torsions. $R'$ is obtained from tubings as in Figure \ref{figure3}, which is isotopic to the result of doing Whitney moves along a clean collection of Whitney disks compatible with $\sx''$(the existence is guaranteed by Fact3)) by Fact 1)2). 

Now Compare the two choices of sheets of $h$ and compute using Proposition \ref{compute}, and we prove the result.\qed

\

\section{proof of the main theorem}\label{sec2}

Recall that an $h$-cobordism is an oriented cobordism $(W;M_0,M_1)$ where the inclusions $M_0\hookrightarrow W$ and $M_1\hookrightarrow W$ are both homotopy equivalences. An $s$-cobordism is an $h$-cobordism with vanishing Whitehead torsion $\tau(W,M_0)$.

\begin{defn}[stabilization of an $h$ cobordism]\label{stab}

\

$W$ is a 5-dimensional $h$-cobordism between $M_0$ and $M_1$. For $\gamma$ a properly smoothly embedded arc from $M_0$ to $M_1$, and a tubular neighborhood $N$ of $\gamma$, we define the stabilization of $W$ along $\gamma$ to be \[(W\backslash N)\cup_{\partial_h N\cong \partial D\times I}(S^2\times S^2-D)\times I\]
where $\partial_h N=\partial N\backslash \operatorname{int}((N\cap(M_0\cup M_1)))$,  and $D$ is a smoothly embedded $4$-dimensional ball in $S^2\times S^2$. $\partial_h N\cong \partial D\times I$ is  induced from a trivialization $N\cong D\times I$.
\end{defn}

Intuitively we want our definition of stabilization to be a connected sum with $S^2\times S^2$ at each level of the $h$-cobordism, however this definition relies on an auxiliary Morse function. We must justify that the diffeomorphism type of the result of a  stabilization is independent of the various choices made in the above definition.
\begin{lem}
Given a 5-dimensional $h$-cobordism $W$ between $M_0$ and $M_1$ and a properly smoothly embedded arc from $M_0$ to $M_1$, the diffeomorphism type of the stabilization of $W$ along $\gamma$ is independent of the choices of $\gamma$, $\gamma$'s tubular neighborhood $N$, smoothly embedded $D\subset S^2\times S^2$ and trivialization $N\cong D\times I$.   
\end{lem}

{\bf Proof:}
We first prove that any two properly embedded arcs $\gamma$ and $\gamma'$ from $M_0$ to $M_1$ are ambient isotopic.

Assume that $\gamma$ and $\gamma'$ have the same boundary points $x_0$ and $x_1$ in $M_0$ and $M_1$ respectively. We have an exact sequence below:
\[\pi_1(M_0\cup M_1,x_0)\to\pi_1(W,x_0)\to\pi_1(W,M_0\cup M_1,x_0)\to\pi_0(M_0\cup M_1,x_0)\]
Since $\pi_1(M_0\cup M_1,x_0)\to\pi_1(W,x_0)$ is an isomorphism,  and $\gamma, \gamma'\in \pi_1(W,M_0\cup M_1,x_0)$ are mapped to the same element in $\pi_0(M_0\cup M_1,x_0)$, $\gamma$ and $\gamma'$ are homotopic with boundary points in $M_0$ and $M_1$ respectively.  

Now in five dimension they are ambient isotopic with boundary points in $M_0$ and $M_1$ respectively by Whitney's embedding lemma.

All smoothly embedded balls in $S^2\times S^2$ are ambient isotopic by the famous theorem of Cerf and Palais.

So we may choose an arbitrary $\gamma$ satisfying the conditions in lemma. By the uniqueness of tubular neighborhood, we may choose an arbitrary tubular neighborhood and a trivialization up to a level-preserving diffeomorphism of $D^4\times I$. It is easily seen that such a diffeomorphism is isotopic to identity(via level-preserving diffeomorphisms). 

Eventually we prove that the different choices made in Definition \ref{stab} give us the same diffeomorphism type of results of stabilizations.\qed

\

Here we need a technical lemma about an $s$-cobordism:
\begin{lem}\label{fund}
For $(W;M_0,M_1)$ a five dimensional $s$-cobordism. For $f$ a compatible Morse function of $W$ with only critical points of index $2,3$, with the images of index $2,3$ critical points to be $\frac{2}{5}$, $\frac{3}{5}$ respectively, we denote $p_1,\hdots,p_k$, $q_1,\hdots, q_k$ its index $2,3$ critical points, $M_\half=f^{-1}(\half)$, and $S_A(p_i), S_D(q_i)\subset M_{\half}$ its ascending spheres and descending spheres. Then the inclusions $M_{\half}\backslash\cup_{i=1}^k S_A(p_i)\hookrightarrow M$, $M_{\half}\backslash\cup_{i=1}^k S_D(q_i)\hookrightarrow M$ both induce isomorphisms on fundamental groups.
\end{lem}
{\bf Proof:}

For critical point $p$ of index $2$, we further denote $D_D(p)$, $D_A(p)$ the descending and ascending disks of $p$ with boundaries in $M_0$ and $M_\half$ respectively.  For critical point $q$ of index $3$, we further denote $D_D(q)$, $D_A(q)$ the descending and ascending disks of $q$ with boundaries in $M_\half$ and $M_1$ respectively.

Note that the inclusion $M_\half\cup_{i=1}^k D_D(q_i)\cup_{i=1}^k D_A(p_i)\hookrightarrow W$ is a homotopy equivalence, so we have $\pi_1(M_\half\cup_{i=1}^k D_D(q_i)\cup_{i=1}^k D_A(p_i))\cong\pi_1(W)$.
By Van Kampen's theorem the inclusion $M_\half\hookrightarrow M_\half\cup_{i=1}^k D_D(q_i)\cup_{i=1}^k D_A(p_i)$ induces an isomorphism on $\pi_1$, since the attached regions are $S^2\times D^2$, which is simply-connected.
So the inclusion $M_\half\hookrightarrow W$ induces isomorphism on $\pi_1$.

Since $S_A(q_i)$ is diffeomorphic to $S^1$, we deduce the inclusions $M_1\backslash\cup_{i=1}^k S_A(q_i)\hookrightarrow M_1\hookrightarrow W$ induce isomorphisms on $\pi_1$.
Note that the natural inclusions $M_\half\backslash\cup_{i=1}^k S_D(q_i)\hookrightarrow W$ and $M_1\backslash\cup_{i=1}^k S_A(q_i)\hookrightarrow W$ differ by an isotopy induced by a gradient-like flow, so the inclusion $M_\half\backslash\cup_{i=1}^k S_D(q_i)\hookrightarrow W$ also induces isomorphism on $\pi_1$.

Now with the above isomorphisms we deduce that the inclusion $M_\half\backslash\cup_{i=1}^k S_D(q_i)\hookrightarrow M_\half$ induces isomorphisms on $\pi_1$. The proof for ascending spheres is the same.
\qed

\

The technical theorem we need is Theorem \ref{tech0}, and to prove which we need a lemma before:

\begin{lem}\label{inj}
$M$ is a four manifold, and $R$ and $S$ are two embedded spheres with trivial normal bundles which intersect transversely in exactly one point $z$. Then the inclusion induces isomorphisms $\pi_1(M\backslash (R\cup S))\to \pi_1(M)$ and $\pi_1(M\backslash R)\to \pi_1(M)$, as well as injections $\pi_2(M\backslash(R\cup S))\to \pi_2(M)$ and $\pi_2(M\backslash R)\to \pi_2(M)$.
\end{lem}

{\bf Proof:}
The injectivity of $\pi_1$ is indeed \cite[Lemma 2.2]{DG}. Assume $\gamma$ is a loop in $M\backslash (R\cup S)$ bounding a immersed disk $D$ in $M$. Assume $D$ intersect $R\cup S$ transversely in some points(all not $z$). We use $R$(or $S$) to tube off the intersections between $D$ and $S$(or $R$), then we get a disk bounding $\gamma$ in $M\backslash(R\cup S)$, and this proves injectivity of $\pi_1$. By general position argument we easily see that $\pi_1(M\backslash (R\cup S))\to \pi_1(M)$ is surjective. The proof of $\pi_1(M\backslash R)\to \pi_1(M)$ is an isomorphism is the same as above.

Denote $\widetilde{M}$ the universal cover of $M$, and $N=\widetilde{M}\backslash(\cup_\alpha (R_\alpha\cup S_\alpha))\subset \widetilde{M}$ the universal cover of $M\backslash (R\cup S)$, where $\alpha$ is indexed by $\pi_1$, and $\#R_\alpha\cap S_\alpha=1$. Denote $N_\alpha$ the tubular neighborhood of $R_\alpha\cup S_\alpha$ in $\widetilde{M}$, then $\partial(N_\alpha)$ is diffeomorphic to $S^3$.
We have exact sequence:
\begin{displaymath}
	H_2(\cup_\alpha {\partial N_\alpha})\to H_2(N)\oplus H_2(N)\to H_2(\widetilde{M})
	\end{displaymath}
	Since $H_2(S^3)=\{0\}$, $H_2(N)\to H_2(\widetilde{M})$ is injective. By Hurwicz's theorem we have that $\pi_2(M\backslash(R\cup S))\to \pi_2(M)$ is an injection.\\
	For $\pi_2(M\backslash R)\to \pi_2(M)$ case we run the exact sequence:
\begin{displaymath}
	H_2(\cup_\alpha {\partial T_\alpha})\to H_2(\cup_\alpha T_\alpha)\oplus H_2(\widetilde{M}\backslash{\cup_\alpha R_\alpha})\to H_2(\widetilde{M})
\end{displaymath}
where $T_\alpha$ is a tubular neighborhood of $R_\alpha$. Then $H_2(\cup_\alpha {\partial T_\alpha})\to H_2(\cup_\alpha T_\alpha)$ is an isomorphism, so $H_2(\widetilde{M}\backslash{\cup_\alpha R_\alpha})\to H_2(\widetilde{M})$ is again an injection, and the injectivity on $\pi_2$ follows.\qed

\begin{thm}\label{tech}
	
	\
	
Given an oriented $4$-dimensional manifold $M$ and smoothly embedded spheres $S_0^i$ and $S_1^i(i=1,\hdots,n)$ satisfying for every $i,j=1,\cdots,n$:
\begin{enumerate}
\item[C1]$\lambda_2(S_0^i,S_1^j)=\pm\delta_{i,j}[1]$, where $1$ represents the trivial element in $\pi_1$. 

\item[C2]$S_0^i$ has trivial normal bundle.

\item[C3]$S_0^{i}\cap S_0^{j}=\varnothing$, $S_1^{i}\cap S_1^{j}=\varnothing$ for $i\ne j$. 

\item[C4]$\pi_1(M\backslash \cup_{i=1}^{n} S_0^i)\to\pi_1(M)$ is injective. 
\end{enumerate}
Then there is another family of spheres $S_3^i(i=1,...,n)$ satisfying for $i,j=1,\cdots,n$:
\begin{enumerate}
\item[R1]$S_3^i(i=1,\hdots,n)$ are disjointly embedded.
	
\item[R2]$S_3^i$ intersects $S_0^j$ transversely and $\#S_3^i\cap S_0^j=\delta_{ij}$. 

\item[R3]$S_3^i$ is homotopic to $S_1^i$.

\item[R4]The Freedman-Quinn invariant $\fq(S_1^i, S_3^i)=0$ for each $i=1,\hdots,n$. 
\end{enumerate}
 Moreover the resulted family $S_3$ is unique up to isotopy, that is if $S_3'^i(i=1,\hdots,n)$ is another family of spheres satisfying the above four conditions, then there exists an isotopy from $\cup_{i=1}^{n}S_3^i$ to $\cup_{i=1}^{n}S_3'^i$.
\end{thm}

{\bf Proof of Theorem \ref{tech}:}

{\bf Step1}: We pick a pair of extra intersections $p,q$ of $S_0^i$, $S_1^j$ for some $i,j$, and $\alpha$, $\beta$ connecting $p$, $q$ in $S_0^i$, $S_1^j$ respectively such that:
\begin{enumerate}
\item[1)]$\alpha$ and $\beta$ do not pass other intersection points of $(\cup_i S_0^i)\cap(\cup_i S_1^i)$ other than their endpoints. 

\item[2)]$\alpha\cup\beta$ is null homotopic in $M$.
\end{enumerate}
The second condition can be satisfied because we can find $p,q$ with opposite intersection number with coefficient in fundamental group given C1.

{\bf Step2}: We find Whitney disk $D$ satisfying:
\begin{enumerate}
\item[1)]$\partial D=\alpha\cup\beta$.

\item[2)]$D$ is transverse to $S_0^i$ and $S_1^i$ on boundary.
\item[3)]int($D$)$\cap (\cup_i S_0^i)=\varnothing$.
\end{enumerate}
The third condition can be achieved because of C4.

{\bf Step3}: Notice that a boundary twist (see, for example, \cite[\S1.3]{FQ} change the framing of a Whitney disk by $2$. Since the algebraic intersections at $p, q$ are of opposite signs, we can perform boundary twist of $D$ on $\beta$ to make the Whitney disks \emph{framed}. The resulted Whitney disks still have interior disjoint from $\cup_i S_0^i$ because $\beta\subset S_1^j$.

{\bf Step4}: We perform Whitney moves to $S_1^i$ along \emph{immersed} $D$ to produce \emph{immersed} spheres $*S_2'^i$, then $\#*S_2'^i\cap(\cup_{j=1}^{n}S_0^j)$ decreases by $2$.

{\bf Step5}: Repeat Step1-Step4 for other extra intersections we obtain immersed spheres $S_2'^i$($i=1,\hdots n$) satisfying:
\begin{enumerate}
	\item[1)]$S_2'^i$ are immersed spheres.
	\item[2)]$\#S_2'^i\cap S_0^j=\delta_{ij}$. 
	\end{enumerate}
{\bf Step5}: For every $j$, we view $S_0^{j}$ as geometric dual of $S_2'^j$, and we use $S_0^{j}$ to tube off:
\begin{enumerate} 
\item[1)] the self-intersection points of $S_2'^j$.

\item[2)] the intersection points of $S_2'^i$ and $S_2'^j$ for $i<j$.
\end{enumerate}
to get disjointly embedded spheres $S_2^i(i=1,\hdots,n)$.

For $j=1,\cdots,n$, we have already known that $S_1^{j}$ and $S_2'^{j}$ are homotopic, to make $S_2'^j$ and $S_2^j$ homotopic, we investigate the effect of tubing the self-intersections of $S_2'^j$ to $S_0^j$.
As we have noted in Proposition \ref{propeff} of Observation \ref{eff}, for a choice of sheets $\sx$ for $S_2'^j$, at a self-intersection $p$, if we tube along an arc as in Proposition \ref{propeff}, then the tubings change the homotopy class of $S_2'^j$ by $\epsilon(p)g_\sx(p)\cdot [S_0^j]$. Now the total change from $S_2'^j$ to $S_2^j$ after we tube off the all self-intersections is:
\begin{equation}\tag{$***$}
\sum_{\textrm{$p$ is a self-intersection of $S_2'^j$}}\epsilon(p)g_\sx(p)\cdot [S_0^j]
\end{equation}
Since $S_2'^j$ is obtained from $S_1^j$ by a regular homotopy, $\mu_2(S_2'^j)=0$, and this means that we can choose the sheets appropriately so that $\sum_{\textrm{$p$ is a self-intersection of $S_2'^j$}}\epsilon(p)g_\sx(p)$ is zero. So we can tube the self-intersections appropriately so that the homotopy class of $S_2'^j$ remains the same. 

For a similar reason we can tube the intersections of $S_2'^j$ with $S_2'^i$ for $i<j$ appropriately(with $S_0^i$ as a geometric dual) so that the homotopy class of $S_2'^j$ remains the same.

Now R1, R2, R3 are already satisfied and we adjust the Freedman-Quinn invariant $\fq(S_1^i, S_2^i)$ for each $i=1\hdots,n$ to achieve R4. Assume for some $i$, $\fq(S_1^i, S_2^i)=[t_1+\hdots+t_s]$, where $t_j\in\tors$. We choose disjointly embedded arcs $\gamma_1,\hdots\gamma_s$ representing $t_1,\hdots t_s$ in $\tors$ with endpoints in $S_2^i$, which are disjoint from $S_2^1\hdots, \hat{S_2^i}, \hdots, S_2^n$ and intersect $S_2^i$ only at endpoints. We modify $S_2^i$ into $S_3^i$ as in Construction \ref{action} by arcs $\gamma_1,\hdots\gamma_s$ with $S_0^i$ as geometric dual, so that $\fq(S_2^i,S_3^i)=[t]$, and $S_3^i$ is supported in a tubular neighborhood of $S_2^i\cup\gamma_1\hdots\cup\gamma_s\cup S_0^i$, thus $S_3^i$ still satisfies R1, R2, R3.

{\bf Proof of uniqueness:}
Without loss of generality, we assume for each $i$, $S_3^i$ coincidence with $S_3'^i$ in a small neighborhood of $z_i$ the intersection of $S_3^i$ and $S_0^i$.
We adopt the procedure of the proof of lightbulb theorem for multiple spheres in \cite[Theorem 10.1]{DG}. If $S_3'^i$ is another family of spheres satisfying R1-R4, then $S_3'^1$ is homotopic to $S_3^1$ in $M$. By Lemma \ref{inj}, if we denote $N=M\backslash(\cup_{i=2}^n S_0^i)$, then $S_3^1$ and $S_3'^1$ are homotopic in $N$. Since we can calculate $\fq(S_3'^1,S_3^1)$ in $N\subset M$, and $\pi_1(N)\to\pi_1(M)$ is an isomorphism(by Lemma \ref{inj} again), we have that $\fq(S_3^1,S_3'^1)=0$ with $N$ as the ambient manifold. Now by lightbulb theorem in non-simply connected case \ref{equal} and the last sentence of Theorem \ref{transitive}, we have that $S_3'^1$ is ambient isotopic to $S_3^1$ in $N\subset M$, and such an isotopy can be chosen to be compactly supported and be supported away from a neighborhood of $S_0^1$. If we view such an isotopy to be an isotopy of $\cup_{i=1}^n S_3'^i$ in $M$, then the isotopy is supported in a compact subset of $M\backslash (\cup_{i=1}^{n}(S_0^i))$. Now we can verify that the new $S_3'^i$ still satisfy R1, R2 with $S_3'^i$ homotopic to $S_3^i$ and $\fq(S_3'^i,S_3^i)=0$, with the new $S_3'^1$ coincident with $S_3^1$. 

If we have arranged that $S_3'^1,\hdots,S_3'^m$ are coincident with $S_3^1,\hdots S_3^m$ respectively and the resulted $S_3'^i(i=1,\hdots n)$ still satisfy R1, R2 with $S_3'^i$ homotopic to $S_3^i$ and $\fq(S_3'^i,S_3^i)=0$, then we apply the argument before to $S_3'^{m+1}$ with $M'=M\backslash (\cup_{i=1}^{m}(S_0^i\cup S_3^i))$. By Lemma \ref{inj} $\pi_1(M')\to\pi_1(M)$ is an isomorphism, and $\pi_2(M')\to\pi_2(M)$ is injective, so both $S_3'^j, S_0^j(j=m+1\hdots,n)$ and $S_3'^j, S_0^j(j=m+1,\hdots n)$ satisfy R1, R2 with $S_3'^i$ homotopic to $S_3^i$ and $\fq(S_3'^i,S_3^i)=0$, \emph{with $M'$ as the ambient manifold}. Now by the same argument for $m=1$ case, we work in the ambient manifold $M'$ with $S_3'^j, S_3^j, S_0^j(j=m+1,\hdots n)$, then we can prove the $m+1$ case and the resulted $n$ spheres still satisfy R1, R2 with $S_3'^j$ homotopic to $S_3^j$ and $\fq(S_3'^j,S_3^j)=0$. Finally we get a new family $S_3'^i(i=1,\hdots,n)$ coincident to $S_3^i(i=1,\hdots,n)$.  \qed 

\

We state the definition of gradient-like vector field following \cite[Definition 3.1]{MIL} here.
\begin{defn}\label{like}
	
	\
	
	Given a cobordism $W^n$ between $M_0^{n-1}$ and $M_1^{n-1}$, and a Morse function $f$ compatible with $W$. Assume further that $W$ is endowed a Riemannian metric. We say that a vector field $\xi$ over $W$ is a \emph{gradient-like vector field} for $f$ if it satisfies the following conditions:
	\begin{enumerate}
		\item[1)]$\xi(f)>0$ holds everywhere except the critical points of $f$.
		\item [2)]Given any critical point $w$ of index $i$, there is a coordinate system in a neighborhood $U$ of $w$ which maps $w$ to $0$: 
		\[(\vec{x},\vec{y})=(x_1,\hdots,x_i,y_1,\hdots,y_{n-i})\]
		so that $f=f(w)-|\vec{x}|^2+|\vec{y}|^2$ and $\xi$ has coordinates $(-x_1,\hdots,-x_i,y_1,\hdots,y_{n-i})$ throughout $U$.
	\end{enumerate}
\end{defn}

We need a technical lemma concerning adjusting Morse function without changing a gradient-like vector field.
\begin{lem}\label{lift}
	
	\

	Given a (Riemannian) cobordism $W^n$ between $M_0^{n-1}$ and $M_1^{n-1}$ and a Morse function $f$ of $W$, and a gradient-like vector field $\xi$ for $f$. For an index $i$ critical point $p$ of $f$, we can find another Morse function $f'(w)=f(w)+\epsilon(w)$, where $\epsilon$ is a non-negative function that can chosen to be supported in any small neighborhood of $p$, so that $f'$ and $f$ has the same collection of critical points, and $\xi$ is still a gradient-like vector field for $f'$.
\end{lem}

{\bf Proof:}
Denote $B(x,r)$ the ball centered at $x\in\mathbb{R}^n$ with radius $r$.
Fix a cut off function $\rho$ on $\mathbb{R}^n$ such that $\rho(x)\in[0,1]$, $\rho(x)=1$ for $x\in B(0,1)$, $|\nabla\rho|\le 1$ and $\rho$ is supported in $B(0,3)$.
Choose a neighborhood $U$ and a coordinate $(\vec{x}, \vec{y})$ as in Definition \ref{like}. Assume $U$ is compact and contains only critical point $p$. For $r$ small enough denote $\delta_r$ the minimum value of $\xi(f)$ on
$f(U)\backslash B(0,r)$. Fix such an $r$, define $\epsilon_r(w)=\frac{1}{2}\delta_r\rho(w/r)$, then $\epsilon_r$ is supported in $B(0,3r)$, and $|\xi(w)(\epsilon_r)|=|\xi(w)\cdot\nabla\epsilon_r(w)|\le|\xi(w)||\nabla\epsilon_r(w)|\le r\cdot(\frac{1}{2}\delta_r/r)=\frac{1}{2}\delta_r$ over $f(U)\backslash B(0,r)$, which means $\xi(f+\epsilon_r)\ge\frac{1}{2}\delta_r$ over $f(U)\backslash B(0,r)$. Set $f_r=f+\epsilon_r$, then $p$ is still the only critical point of $f_r$ on $U$, and $\xi$ is still a gradient-like vector field on of $f_r$.\qed

\

{\bf Proof of The Main Theorem \ref{main}:}

We give $W$ a Riemannian metric. For Morse function $f$ as in Theorem \ref{main}, We further assume that for a critical point $a$ of index $i$, $f(a)=\frac{i}{5}$. Denote the index-$2$ points of $f$ by $p_1,\hdots,p_k$, index-$3$ points of $f$ by $q_1,\hdots,q_k$. Denote $M_{c}=f^{-1}(c)$.

For critical point $p$ of index 2, denote $D_D(p)$, $D_A(p)$ the descending and ascending disks of $p$ with boundaries in $M_0$ and $M_\half$ respectively. Denote $S_D(p)$ and $S_A(p)$ the intersection of $D_D(p)$ with $M_0$ and $D_A(p)$ with $M_\half$ respectively. For critical point $q$ of index 3, denote $D_D(q)$, $D_A(q)$ the descending and ascending disks of $q$ with boundaries in $M_\half$ and $M_1$ respectively. Denote $S_D(q)$ and $S_A(q)$ the intersection of $D_D(q)$ with $M_\half$ and $D_A(q)$ with $M_1$ respectively.    

\

{\bf Step 1:}  

By Proposition \ref{White}, since the Whitehead torsion $\tau(W,M_0)=0$, we adjust $f$ and $\xi$ so that $\lambda_2(S_A(p_i),S_D(q_j))=\delta_{ij}[1]$.

\

{\bf Step2:} 

Now we verify that $M_\half$, $S_D(q_1)$, $S_A(p_1)$ satisfy the conditions in Lemma \ref{tech}.

After {\bf Step 1} condition 1) of Lemma \ref{tech} is already satisfied.

Condition 2) and 3) are properties of a general cobordism.

By Lemma \ref{fund} condition 4) is also satisfied.

Now we can apply Theorem \ref{tech} to $M_\half$, that is we can choose a different family of spheres $S'_D(q_i)$ satisfying for $i,j=1,\hdots,k$:
\begin{enumerate}
\item[1)]$S'_D(q_i)$ are disjointly embedded.
\item[2)]$\#S'_D(q_i)\cap S_A(p_j)=\delta_{ij}$ and the intersections are transverse. 
\item[3)]$S_D(q_i)$ and $S'_D(q_i)$ are homotopic.
\item[4)]$\fq(S_D(q_i),S'_D(q_i))=0$.
\end{enumerate}

{\bf Step3}:

Choose a pair of geometric dual spheres in $M_\half$, say $S'_D(q_1)$ and $S_A(p_1)$. We work in the following substeps to turn $S_A(p_1)$ into a geometric dual of \emph{both} $S_D(q_1)$ and $S_D'(q_1)$ \emph{in the stabilized middle level}. The proof is essentially part of the proof in \cite{AD}.

After {\bf Step 1}, as is the proof of Lemma \ref{tech}, we pick out a canceling intersection pair(if any) $u,v$ of $S_D(q_1)$ and $S_A(p_1)$, so that there are Whitney arcs $\alpha$ and $\beta$ connecting $u$ and $v$ in $S_D(q_1)$ and $S_A(p_1)$ respectively, and we can find immersed Whitney disks $D$ with boundary $\alpha\cup\beta$, and $\operatorname{int}(D)$ are disjoint from $S_D(q_1)$. We perform boundary twists(see, for example \cite[\S 1.3]{FQ}) at $\beta$ so $D$ are all framed.

\begin{figure}[!htb]
\centering
\subfigure{\resizebox*{12cm}{!}{\includegraphics{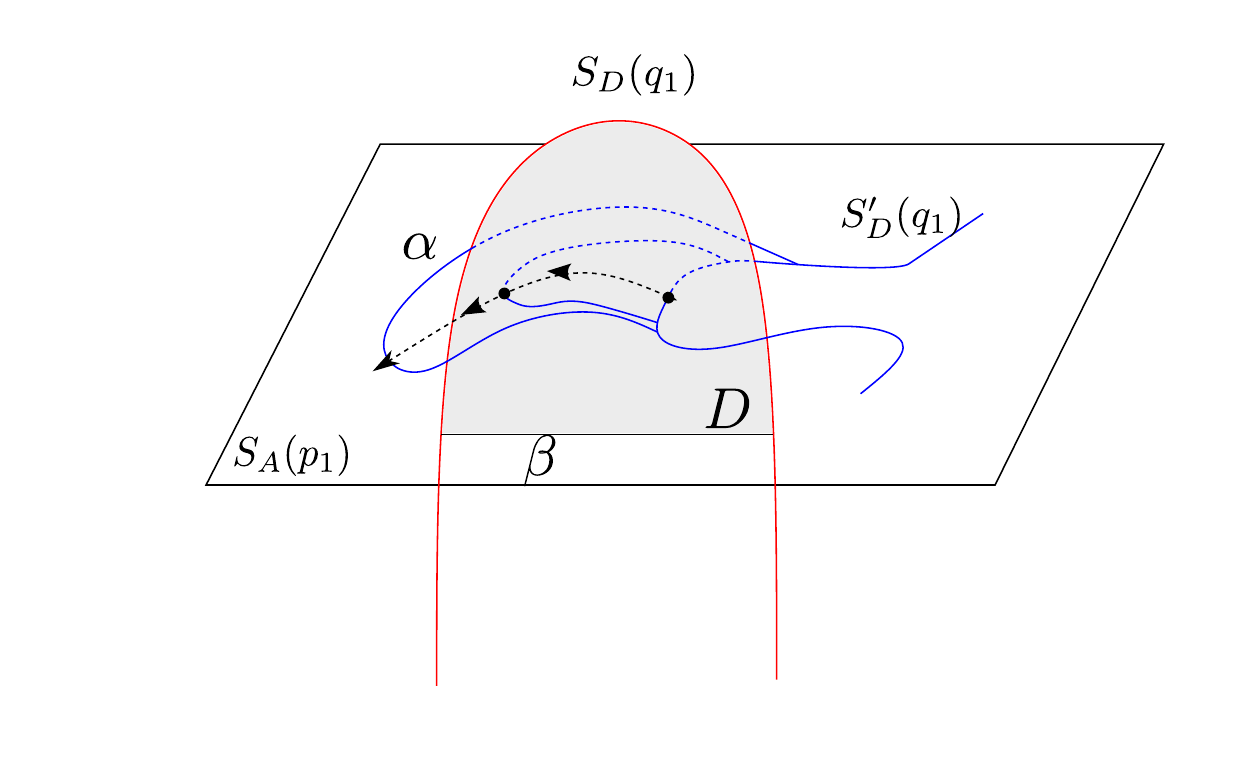}}}
\caption{}\label{figure4}
\end{figure}

We want to perform Whitney moves on $S_A(p_1)$ to cancel the extra intersections $u,v$ of $S_A(p_1)$ and $S_D(q_1)$, however this may create new intersections between $S_A(p_1)$ and $S'_D(q_1)$. To mend this, for each intersection between $D$ and $S'_D(q_1)$, we perform finger moves on $S'_D(q_1)$ across $D$ towards $\alpha$ to obtain embedded $*S''_{D}(q_1)$. The finger move is indicated in Figure \ref{figure4}. 
We must do the finger moves appropriately:
\begin{equation}\tag{$****$}\label{appropriate}
		\begin{split}
&\emph{All guiding curves of the finger moves are mutually disjoint, which are disjoint from $S_A(p_1)$}\\
&\emph{and intersects $*S''_{D}(q_1)$$($or $S'_D(q_1)$ when we are canceling the first pair of $u,v$$)$ in a single}\\
&\emph{endpoint}.
\end{split}
\end{equation}

Then $*S''_{D}(q_1)$ is isotopic to $S'_D(q_1)$ and is disjoint from $D$ and still intersects $S_A(p_1)$ transversely in exactly one point.

We perform Whitney moves on $S_A(p_1)$  along $D$ to obtain an \emph{immersed} sphere $*S$ so that $*S$ intersects $*S''_{D}(q_1)$ transversely in exactly one point and intersects $S_D(q_1)$ in two points less. Repeat the above procedures we find an immersed sphere $S$, as well as $S''_D(q_1)$ which is isotopic to $S'_D(q_1)$, and $S$ intersects both $S_D(q_1)$ and $S''_D(q_1)$ in exactly one point, and by \ref{appropriate}, $S''_D(q_1)$ intersects $S_A(p_1)$ in exactly one point.

From now on we work in the stabilized cobordism $W'$(along an appropriately chosen arc $\gamma$, see Definition \ref{stab}) with extended Morse function and gradient-like vector field,  which we still denote by $f$ and $\xi$. We denoted the new middle level manifold by $M'_\half=M_\half\#\stb$. We tube $S$ along an arc to the first component of $S^2\times S^2$, say $S^2\times\{a\}$ to obtain $S'$. We can use $\{b\}\times S^2$ as a geometrical dual to $S'$ to tube off the self-intersections of $S'$ to obtain an embedded $S''$. Notice that by Proposition \ref{propeff}, the homology class of $S''$ in $M'_\half$ can be calculated as (if we orient $\stb$, $S^2\times\{a\}$ and $\{b\}\times S^2$ so that $[S^2\times\{a\}]\cdot[\{b\}\times S^2]=1$) $[S'']=[S]+[S^2\times\{a\}]+n[\{b\}\times S^2]$
, where $n$ is the self intersection number of $S$(or $S'$) counted with sign. Since $S$ is obtained from $S_A(p_1)$ by regular homotopies, we have $[S]^2=[S_A(p_1)]^2=0$ and $n=0$. With above equalities we calculate that $[S'']^2=0$, so $S''$ has trivial normal bundle.

Now $S_D(q_1)$ and $S''_D(q_1)$ has a common geometric dual $S''$ in $M'_\half$ , and as  a result of Lemma \ref{tech}, $S_D(q_1)$, $S'_D(q_1)$ and $S''_D(q_1)$ are mutually homotopic.

\

{\bf Step 5}: 
Now $0=\fq(S_D(q_1),S'_D(q_1))=\fq(S_D(q_1),S''_D(q_1))$(the first equality holds because of Step 2 and 4) of Remark \ref{remofFQ}(the second equality holds because $S'_D(q_1)$ is isotopic to $S''_D(q_1)$ by Step 3. By Theorem \ref{equal} $S_D(q_1)$ and $S''_D(q_1)$ are isotopic in $M'_\half$.
\

{\bf Step 6:} Cancel a pair of critical points in the once-stabilized cobordism.

Now $S_D(q_1)$ is isotopic to $S''_D(q_1)$ is $M'_\half$, by \cite[Lemma 4.6]{MIL}, we can alter $\xi$ by $\xi'$, which is still a gradient-like vector field of $f$, such that the new descending sphere of $q_1$ is exactly $S''_D(q_1)$ and the ascending sphere of $p_1$ is still $S_A(p_1)$. 

Now by Lemma \ref{lift}, we can increase $f'$ in a small neighborhood of $p_1$ and decrease $f'$ in a small neighborhood of $q_1$ to get $f''$, so that $\xi'$ is still a gradient-like vector field for $f''$. Now $f''$ decomposes $W'$ into a composition of cobordisms $c_+\circ c\circ c_-$, where $c$ is the middle cobordism with only critical points $p_1$, $q_1$, and their ascending sphere and descending sphere have only one intersection. Now by the cancellation theorem of Milnor \cite[Theorem 5.4]{MIL}, we can adjust $f''$ and $\xi'$ so that by $g$ and $\eta$ so that $g$ has $k-1$ pairs of critical points, and $\eta$ is a gradient-like vector field for $g$. 

Apply the above procedures to other pairs of critical points, and finally we get a trivial cobordism after $k$ times of stabilizations.\qed

\

\section{Remarks of the Main Theorem}\label{app}
		In this section  we simplify the proof of the main theorem in the case where $W$ is simply-connected, however this approach will only prove that after $k+1$ times of stabilizations we can get a product cobordism.
		In this section we will use a different notation of cobordism:
		\begin{defn}\label{cob}
			
			\
			
			A \emph{five dimensional cobordism between oriented four manifolds} $M_0$ and $M_1$ consists of a five dimensional manifold $W$, and inclusions $i_0:M_0\to W$ and $i_1:M_1\to W$ with $\partial W=i_0(M_0)\sqcup \overline{i_1(M_1)}$, which we will denote by  $(W;M_0,M_1;i_0,i_1)$. An \emph{equivalence} between $h$-cobordisms $(W^1;M_0^1,M_1^1;i_0^1,i_0^2)$ and $(W^2;M_0^2,M_1^2;i_0^2,i_1^2)$ is a  triple $(\phi,\phi_1,\phi_2)$,where $\phi:W^1\to W^2 $, $\phi_0:M_0^1\to M_0^2$, $\phi_1:M_1^1\to M_1^2$ are all diffeomorphisms, and they commute with $i_0^1,i_1^1,i_0^2,i_1^2$ in the obvious way.   
		\end{defn}
		
		\begin{rem}\label{induce}
			
			\

			Notice that an $h$-cobordism $(W;M_0,M_1;i_0,i_1)$ between $M_0$ and $M_1$ induces an isometry between $Q(M_0)$ and $(Q(M_1))$ via $i_0^*:H_2^*(W)\to H_2^*(M_0)$ and $i_1^*:H_2^*(W)\to H_2^*(M_1)$.
			\end{rem}
			
			\

		We state two lemmas needed in this section. The following lemma from \cite{WT2} is well-known.
		\begin{thm}\label{diff}
			Let $N$ be a simply-connected, closed, smooth four manifold. Suppose that either
			\begin{enumerate}
				\item[1)]$Q(N)$ is indefinite, or
				\item[2)]$\operatorname{rank}(H^2(N))\le8$.
			\end{enumerate}
			Then every isometry of $Q(N\#S^2\times S^2)$ is realized by a diffeomorphism of $N\#\stb$.
		\end{thm}
		
		We also need the classification theorem of simply-connected $h$-cobordisms in \cite{KR}.
		\begin{thm}\label{classify}
			Let $M_0$ and $M_1$ be fixed closed oriented smooth simply connected  four manifolds. Then the set $\{$$h$-cobordisms between $M_0$ and $M_1$$\}/$equivalences of the form $(\phi,\id,\id)$ is isomorphic to the set of morphisms between the intersection forms of $Q(M_0)$ and $Q(M_1)$.
		\end{thm}
		
		\

		{\bf An alternative proof of the Main Theorem \ref{main} in the simply-connected case:}
		
		Assume $W$ is non-trivial, so $k\ge 1$.
		Since $W$ is an $h$-cobordism between $M_0$ and $M_1$ with $k$ pairs of critical points of index $2$,$3$, it is well-known that $M_0\#^k\stb$ is diffeomorphic to $M_1\#^k\stb$. After $k$ times of stabilizations of $W$, we get an $h$-cobordism between $M_0\#^k\stb$ and $M_1\#^k\stb$. Denote the new $h$-cobordism by $(W';M_0',M_1';i_0',i_1')$ under the notation of Definition \ref{cob}, where $M_0'$ and $M_1'$ are regarded as submanifolds of $W'$, and $i_0$ and $i_1$ are inclusions.
		
		Denote the one time stabilization of $(W';M_0',M_1';i_0,i_1)$ by $(W'';M_0'',M_1'';i_0'',i_1'')$. We prove that $W''$ is a product cobordism.
		
		Since $k\ge 1$, $M_0'$ and $M_1'$ both has indefinite intersection forms. By Remark \ref{induce}, $(W'';M_0'',M_1'';i_0'',i_1'')$ induces an isometry $\phi$ between $H^2(M_0'')$ and $H^2(M_1'')$. Now we can pick a diffeomorphism $f:M_0''\to M_1''$ realizing $\phi$ by Theorem \ref{diff}. Now consider the $h$-cobordism $(M_1''\times I;M_0'', M_1''; f\times \{0\},\Id \times\{0\})$, which induces exactly $\phi$ between $H^2(M_0'')$ and $H^2(M_1'')$. By Theorem \ref{classify}, $W''$ is diffeomorphic to $M_1''\times I$. By construction $W''$ is the $k+1$ times stabilizations of $W$, and this ends the proof.\qed
		
		\begin{rem}
			
			\
			
	John W. Morgan and Zolt\'{a}n Szab\'{o} pointed out in \cite{MO} that the family of $h$-cobordisms they constructed in \cite[Theorem 1.1]{MO} has $C'=1$(see Section \ref{sec1} for definition of $C'$), so they can be stabilized to product cobordism after one stabilization by Theorem \ref{main}. Actually any non-trivial $h$-cobordism between two simply connected, closed, \emph{diffeomorphic} four manifolds \emph{with indefinite intersection forms or $\operatorname{rank}(H^2)\le8$} has $C'=1$, which can be proved by the same technique used above.
	\end{rem}
\

\end{document}